\documentclass[12pt]{article}
\usepackage{graphicx}
\usepackage{amsfonts}
\usepackage{amsmath}
\usepackage{amssymb,amsbsy}
\usepackage{latexsym}
\usepackage{epsfig}
\usepackage[textwidth=18cm]{geometry}
\usepackage{nccmath}
\usepackage{revsymb}
\usepackage{eurosym}
\usepackage{latexsym}

\input{tcilatex}
\begin{document}

\title{\textbf{Solution of constrained mechanical multibody systems using
Adomian decomposition method }}
\author{Brahim Benhammouda \thanks{%
Email: bbenhammouda@hct.ac.ae} \\
{\small {Higher Colleges of Technology. Abu Dhabi Men's College} }\\
{\small {P.O. Box 25035}}, {\small {Abu Dhabi, United Arab Emirates.} }\\
}
\maketitle

\textbf{Abstract. }Constrained mechanical multibody systems arise in many
important applications like robotics, vehicle and machinery dynamics and
biomechanics of locomotion of humans. These systems are described by the
Euler-Lagrange equations which are index-three differential-algebraic
equations (DAEs) and hence difficult to treat numerically.

The purpose of this paper is to propose a novel technique to solve the
Euler-Lagrange equations efficiently. This technique applies the Adomian
decomposition method (ADM) directly to these equations. The great advantage
of our technique is that it neither applies complex transformations to the
equations nor uses index-reductions to obtain the solution. Furthermore, it
requires solving only linear algebraic systems with a constant nonsingular
coefficient matrix at each iteration. The technique developed leads to a
simple general algorithm that can be programmed in Maple or Mathematica to
simulate real application problems. To illustrate the effectiveness of the
proposed technique and its advantages, we apply it to solve an example of
the Euler-Lagrange equations that describes a two-link\ planar robotic
system.



\section{Introduction}

Constrained mechanical multibody systems arise in many areas of applications
such as robotics, biomechanics of locomotion of humans and dynamics of
vehicle and machinery \cite{Haug,Truck-model,Knee,Wash-m-1}. The dynamical
behavior of constrained multibody systems is described by the Euler-Lagrange
equations%
\begin{equation}
\left\{ 
\begin{array}{l}
\dfrac{dp}{dt}=v, \\ 
M(p)\dfrac{dv}{dt}=f\left( p,v\right) -G^{^{\text{\textsf{T}}}}\left(
p\right) \lambda \text{,} \\ 
0=g\left( p\right) \text{, \ }t\geq 0\text{.}%
\end{array}%
\right.  \label{M1}
\end{equation}%
Here $t$ is the time, $p\left( t\right) \in \mathbb{R}^{n_{p}}$ and $v\left(
t\right) \in \mathbb{R}^{n_{p}}$ specify the positions and the orientations
of all bodies and their velocities, respectively. The vector $\lambda \in 
\mathbb{R}^{n_{\lambda }}$ is the vector of Lagrange multipliers. The matrix 
$M\left( p\right) \in \mathbb{R}^{n_{p}}\times \mathbb{R}^{n_{p}}$ is a mass
matrix. The mapping $f:\mathbb{R}^{n_{p}}\times \mathbb{R}%
^{n_{p}}\longrightarrow \mathbb{R}^{n_{p}}$\ defines the applied and
internal forces (other than the constraint forces), whereas $g:\mathbb{R}%
^{n_{p}}\longrightarrow \mathbb{R}^{n_{\lambda }},$ ($n_{\lambda }\leq n_{p}$%
) defines the constraints. The term $G(p)^{^{\text{\textsf{T}}}}\lambda $
represents the constraint forces, where $G\left( p\right) =\partial g/dp\in 
\mathbb{R}^{n_{\lambda }\times n_{p}}$ denotes the Jacobian of $g\left(
p\right) $. \ 

The Euler-Lagrange equations (\ref{M1}) form a nonlinear system of
differential-algebraic equations (DAEs). Consistent initial conditions 
\begin{equation}
p\left( 0\right) =p_{0}\text{, \ }v\left( 0\right) =v_{0}\text{,}  \label{M2}
\end{equation}%
are necessary to uniquely determine a solution. The given vectors $p_{0}\ $%
and $v_{0}$, which specify the initial configuration and initial velocity,
are chosen so that the consistency equations 
\begin{eqnarray}
g\left( p_{0}\right) &=&0,  \label{Consistency-1} \\
G(p_{0})v_{0} &=&0,  \label{Consistency-2}
\end{eqnarray}%
are satisfied. For the variable $\lambda \left( t\right) ,$ no initial
condition is prescribed because $\lambda \left( 0\right) $ is already
determined by DAE (\ref{M1}) and initial conditions (\ref{M2}).

Throughout this paper, we assume that $M\left( p\right) $, $f\left(
p,v\right) $ and $g\left( p\right) $ are analytical. A standard assumption
on the Jacobian $G\left( p\right) $ is the full row rank condition%
\begin{equation}
\text{rank }G\left( p\right) =n_{\lambda },  \label{Full-rank}
\end{equation}%
which means that the constraint equations defined by $g(p)$ are linearly
independent.

In addition the matrix $M$ is assumed to be symmetric and positive definite,
that is 
\begin{equation}
z^{^{\text{\textsf{T}}}}M\left( p\right) z>0,\text{for all }z\in \text{Ker }%
G(p).  \label{Positive}
\end{equation}

\noindent If assumptions (\ref{Full-rank})-(\ref{Positive}) hold, then the
matrix 
\begin{equation}
\left( 
\begin{array}{cc}
M\left( p\right) & G^{^{\text{\textsf{T}}}}\left( p\right) \\ 
G\left( p\right) & 0%
\end{array}%
\right) ,  \label{nonsingular}
\end{equation}%
is nonsingular and DAE (\ref{M1}) is therefore index-three. For any
consistent initial conditions (\ref{M2}), system (\ref{M1}) has a unique
solution \cite{Haug}. The index we consider here is the differential index
which is the minimum number of times that all or a part of the DAE must be
differentiated with respect to time in order to obtain an ordinary
differential equation \cite{Index}. It is well-known that index-three DAEs
present difficulties for numerical integration methods \cite{Brenan}.
Therefore, several techniques have been proposed in the literature to solve
DAE system (\ref{M1}) \cite{Brenan,Baumgarte,Stab-0,Stab-1,
Robot-example,Brahim-SC}. A very popular way is to reduce the index by
differentiating the constraints one or more times with respect to time
before applying numerical integration methods. However, the main problem
with this technique is that the numerical solution of the index-reduced
system may no longer satisfy the constraints of original DAE (\ref{M1}) due
to error propagation. Constraints violation, known also as drift-off
phenomena, leads to non physical solutions. To overcome this difficulty,
some techniques like stabilization or augmented Lagrangian formulation have
been proposed to keep the constraint violations under control during the
numerical integration \cite{Baumgarte,Stab-0,Stab-1,Stab-2}. The most
popular stabilization method is that of Baumgarte \cite{Baumgarte}, but its
drawback is the way of choosing its feedback parameters. The augmented
Lagrangian formulation \cite{Aug-Lag} has the same problem of parameter
selection. The challenge is therefore to construct efficient methods that
provide solutions of the Euler-Lagrange equations which satisfy the
constraints in these equations.

The Adomian decomposition method (ADM) and its modifications \cite{Alge-ADM,
ADM-ode-1, ADM-ode-2, ADM-ode-3, ADM-ode-4, ADM-Mod-1, ADM-ode-Mod-2} are
known to be efficient methods in solving a large variety of linear and
nonlinear problems in science and engineering. Among these problems, we
mention algebraic equations \cite{Alge-ADM}, ordinary differential equations 
\cite{ADM-ode-1, ADM-ode-2, ADM-ode-3, ADM-ode-4, ADM-Mod-1, ADM-ode-Mod-2},
partial differential equations \cite{ADM-PDEs} and integral equations \cite%
{ADM-Mod-Integrals}.

In this work, we present a new approach to solve the Euler-Lagrange
equations using ADM. The solution by this method satisfies all the DAE
constraints. The ADM\ is first applied directly to the Euler-Lagrange
equations where the nonlinear terms are expanded using the Adomian
polynomials \cite{ADP-1,ADP-2,ADP-3,ADP-4,ADP-5,ADP-6}. Based on the index
of the Euler-Lagrange equations, a nonsingular linear algebraic recursion
system is derived for the expansion components of the solution. Our
technique has the great advantage that it does not use complex
transformations like index reductions before applying the ADM to the
equations. To demonstrate the effectiveness of the proposed method, we solve
an example of the Euler-Lagrange equations that models a two-link planar
robotic system. Further, our technique is based on a simple algorithm that
can be programmed in Maple or Mathematica to simulate real application
problems.

This paper is organized as follows: in section \ref{ADM}, we review the ADM
for solving ordinary differential equations. Next, in section \ref{Multibody
system-ADM}, we present our method for the solution of the Euler-Lagrange
equations. Then, in section \ref{Test problems}, we apply the developed
technique to solve an example of the Euler-Lagrange equations that models a
two-link \ planar robotic system. Finally, a discussion and a conclusion are
given in sections \ref{Discussion} and \ref{Conclusion}, respectively.

\section{Adomian decomposition method}

\label{ADM}In this section, we give a brief review for the Adomian
decomposition method (ADM) \cite{Alge-ADM, ADM-ode-1, ADM-ode-2, ADM-ode-3,
ADM-ode-4, ADM-Mod-1, ADM-ode-Mod-2} to solve ordinary differential
equations. For this purpose, let us consider the following nonlinear
differential equation%
\begin{equation}
Lu+Ru+N(u)=f,  \label{ADM-1}
\end{equation}%
where $L$ is an easily invertible operator (usually taken as the
highest-order derivative), $R$ is an operator grouping the remaining
lower-order derivatives, $N\left( u\right) $ is the nonlinear term and $f$
is a given analytical function.

Solving equation (\ref{ADM-1}) for $Lu$ then applying the inverse operator $%
L^{-1}$ to both sides, we obtain%
\begin{equation}
L^{-1}Lu=L^{-1}f-L^{-1}Ru-L^{-1}N\left( u\right) .  \label{ADM-2}
\end{equation}%
If $Lu=du/dt$ and the initial condition $u\left( t_{0}\right) =u_{0}$ is
given, then $L^{-1}$ represents the integral from $t_{0}$ to $t$ and $%
L^{-1}Lu=u-u_{0}.$

\begin{equation}
u=u_{0}+L^{-1}f-L^{-1}Ru-L^{-1}N\left( u\right) .  \label{ADM-3}
\end{equation}%
To apply the ADM to equation (\ref{ADM-3}), we first assume that the
solution $u$ of (\ref{ADM-1}) to have the infinite series form 
\begin{equation}
u=\sum_{n=0}^{\infty }u^{\left( n\right) },  \label{ADM-4}
\end{equation}%
where the unknown solution components $u^{\left( n\right) },$ $%
n=0,1,2,\ldots $ are to be determined later by the method.

Second, the nonlinear term $N(u)$ is expanded in an infinite series in terms
of the Adomian polynomials $N^{\left( n\right) }$ \cite%
{ADP-1,ADP-2,ADP-3,ADP-4,ADP-5,ADP-6} as 
\begin{equation}
N\left( u\right) =\sum_{n=0}^{\infty }N^{\left( n\right) }\left( u^{\left(
0\right) },\ldots ,u^{\left( n\right) }\right) .\text{ }  \label{ADM-5}
\end{equation}%
Substituting (\ref{ADM-4}) and (\ref{ADM-5}) into (\ref{ADM-3}) and choosing 
$u^{\left( 0\right) }$ as 
\begin{equation}
u^{\left( 0\right) }=u_{0}+L^{-1}f\text{,}  \label{ADM-6}
\end{equation}%
we obtain%
\begin{equation}
\sum_{n=0}^{\infty }u^{\left( n\right) }=u^{\left( 0\right)
}-L^{-1}R\sum_{n=0}^{\infty }u^{\left( n\right) }-L^{-1}\sum_{n=0}^{\infty
}N^{\left( n\right) }.  \label{ADM-7}
\end{equation}%
Comparing the general term on the left hand side with that on the right hand
side, we derive the following recursion scheme for the ADM 
\begin{equation}
u^{\left( n\right) }=-L^{-1}Ru^{\left( n-1\right) }-L^{-1}N^{\left(
n-1\right) },\text{ \ }n\geq 1.  \label{ADM-8}
\end{equation}%
Since $u^{\left( 0\right) }$ is known, recursion (\ref{ADM-8}) can be used
to generate as many solution components $u^{\left( n\right) }$ as one wants.
Further, if series (\ref{ADM-4}) converges then it gives the exact solution
of (\ref{ADM-1}) and an approximation of order $n_{0}$ to solution can be
obtained from 
\begin{equation}
u=\sum_{n=0}^{n_{0}-1}u^{\left( n\right) }.  \label{ADM-9}
\end{equation}%
To compute the Adomian polynomials $N^{\left( n\right) }$, $n=0,1,\ldots $
associated with the nonlinearity $N\left( u\right) ,$ one can use the
following definition for all forms of nonlinearity 
\begin{equation}
N^{\left( n\right) }:=N^{\left( n\right) }\left( u^{\left( 0\right) },\ldots
,u^{\left( n\right) }\right) =\frac{1}{n!}\frac{d^{n}}{d\lambda ^{n}}\left(
N\left( \sum_{i=0}^{\infty }\lambda ^{i}u^{\left( i\right) }\right) \right)
_{\lambda =0}\text{, \ }n\geq 0.  \label{ADM-10}
\end{equation}%
Using this formula, we obtain the following first few Adomian polynomials 
\begin{eqnarray}
N^{\left( 0\right) } &=&N\left( u^{\left( 0\right) }\right) ,  \notag \\
N^{\left( 1\right) } &=&u^{\left( 1\right) }N^{\prime }\left( u^{\left(
0\right) }\right) ,  \notag \\
N^{\left( 2\right) } &=&u^{\left( 2\right) }N^{\prime }\left( u^{\left(
0\right) }\right) +\dfrac{\left( u^{\left( 1\right) }\right) ^{2}}{2!}%
N^{\prime \prime }\left( u^{\left( 0\right) }\right) ,  \label{ADM-11} \\
N^{\left( 3\right) } &=&u^{\left( 3\right) }N^{\prime }\left( u^{\left(
0\right) }\right) +u^{\left( 1\right) }u^{\left( 2\right) }N^{\prime \prime
}\left( u^{\left( 0\right) }\right) +\dfrac{\left( u^{\left( 1\right)
}\right) ^{3}}{3!}N^{\prime \prime \prime }\left( u^{\left( 0\right)
}\right) ,  \notag \\
N^{\left( 4\right) } &=&u^{\left( 4\right) }N^{\prime }\left( u^{\left(
0\right) }\right) +\left( \dfrac{\left( u^{\left( 2\right) }\right) ^{2}}{2!}%
+u^{\left( 1\right) }u^{\left( 3\right) }\right) N^{\prime \prime }\left(
u^{\left( 0\right) }\right) +\dfrac{\left( u^{\left( 1\right) }\right)
^{2}u^{\left( 2\right) }}{2!}N^{\prime \prime \prime }\left( u^{\left(
0\right) }\right) +\dfrac{\left( u^{\left( 1\right) }\right) ^{4}}{4!}%
N^{\prime \prime \prime \prime }\left( u^{\left( 0\right) }\right) ,  \notag
\end{eqnarray}%
where the dash $(^{\prime })$ represents the differentiation with respect to 
$u$.

In a similar manner, one can easily generate the remaining polynomials from (%
\ref{ADM-10}). In the literature, there are several algorithms for computing
the Adomian polynomials without the need for formula (\ref{ADM-10}), but a
more convenient algorithm for the $m$-variable case is recently proposed in 
\cite{ADP-6} 
\begin{equation}
N^{\left( n\right) }=\dfrac{1}{n}\sum_{k=1}^{m}\sum_{i=0}^{n-1}\left(
i+1\right) v_{k}^{\left( i+1\right) }\frac{\partial N^{\left( n-1-i\right) }%
}{\partial v_{k}^{\left( 0\right) }},\text{ \ }n\geq 1.  \label{ADM-12}
\end{equation}

\section{The proposed method}

\label{Multibody system-ADM}In this section, we present our method for
solving the Euler-Lagrange equations. These equations are known to be
difficult to treat numerically since they represent an index-three system of
differential-algebraic equations (DAEs). The technique we propose here is
based on the Adomian decomposition method (ADM). To solve these equations,
we first apply the ADM directly to them and expand the nonlinear terms using
the Adomian polynomials. Then, an algebraic linear recursion system for the
solution expansion components is derived. Taking account of the index of the
DAE, this algebraic system is shown to be uniquely solvable for the solution
expansion components. The main advantage of our technique is that it does
not require to transform the equations to lower index DAEs before applying
the ADM to them. We start our technique by the following proposition.

\textbf{Proposition}: Let $u=\left( u_{1},\ldots ,u_{l}\right) ^{^{\text{%
\textsf{T}}}}\in \mathbb{R}^{l},$ $v=\left( v_{1},\ldots ,v_{n}\right) ^{^{%
\text{\textsf{T}}}}$ $\in \mathbb{R}^{n}$ be two vectors, $A:=(A\left(
u\right) _{i,j})\in \mathbb{R}^{m}$ $\times \mathbb{R}^{n}$ be a matrix, $%
u=\sum_{k=0}^{\infty }u^{\left( k\right) },v=\sum_{k=0}^{\infty }v^{\left(
k\right) },$ $u^{\left( k\right) }=$ $\left( u_{1}^{\left( k\right) },\ldots
,u_{l}^{\left( k\right) }\right) ,$ $v^{\left( k\right) }=$ $\left(
v_{1}^{\left( k\right) },\ldots ,v_{n}^{\left( k\right) }\right) .$ Let $%
A^{\left( k\right) }$ $=\left( A_{i,j}^{\left( k\right) }\right) ,$ where $%
A_{i,j}^{\left( k\right) }=A_{i,j}^{\left( k\right) }\left( u^{\left(
0\right) },\ldots ,u^{\left( k\right) }\right) $ denotes the $k$-th Adomian
polynomial of the entry $A_{i,j}$. Then, given $A$ and $v$, the Adomian
polynomials of $z=Av$ and $z=A^{-1}v$ are given by

(a) $z^{\left( k\right) }=\sum_{l=0}^{k}A^{\left( l\right) }v^{\left(
k-l\right) }=$ $\sum_{l=0}^{k}A^{\left( k-l\right) }v^{\left( l\right) }$ and

(b) $A^{\left( 0\right) }z^{\left( k\right) }=v^{\left( k\right)
}-\sum_{l=0}^{k-1}A^{\left( k-l\right) }z^{\left( l\right) },$ $k=0,1,\ldots 
$

\textbf{Proof}:

(a) Let $z=Av$, then the Adomian polynomials $z^{\left( k\right) }$ of $z$
are given by 
\begin{eqnarray*}
z^{\left( k\right) } &=&\left( Av\right) ^{\left( k\right) }=\left( 
\begin{array}{c}
\left( \sum_{j=1}^{n}A_{1,j}v_{j}\right) ^{\left( k\right) } \\ 
\vdots  \\ 
\left( \sum_{j=1}^{n}A_{m,j}v_{j}\right) ^{\left( k\right) }%
\end{array}%
\right) =\left( 
\begin{array}{c}
\sum_{j=1}^{n}\left( A_{1,j}v_{j}\right) ^{\left( k\right) } \\ 
\vdots  \\ 
\sum_{j=1}^{n}\left( A_{m,j}v_{j}\right) ^{\left( k\right) }%
\end{array}%
\right) , \\
&=&\left( 
\begin{array}{c}
\sum_{j=1}^{n}\sum_{l=0}^{k}A_{1,j}^{\left( k-l\right) }v_{j}^{\left(
l\right) } \\ 
\vdots  \\ 
\sum_{j=1}^{n}\sum_{l=0}^{k}A_{m,j}^{\left( k-l\right) }v_{j}^{\left(
l\right) }%
\end{array}%
\right) =\left( 
\begin{array}{c}
\sum_{l=0}^{k}\sum_{j=1}^{n}A_{1,j}^{\left( k-l\right) }v_{j}^{\left(
l\right) } \\ 
\vdots  \\ 
\sum_{l=0}^{k}\sum_{j=1}^{n}A_{m,j}^{\left( k-l\right) }v_{j}^{\left(
l\right) }%
\end{array}%
\right) , \\
&=&\sum_{l=0}^{k}A^{\left( k-l\right) }v^{\left( l\right) }.
\end{eqnarray*}

(b) Let $Az=v$ then, using (a), the Adomian polynomials $v^{\left( k\right) }
$ of $v$ are given by 
\begin{equation*}
v^{\left( k\right) }=\left( Az\right) ^{\left( k\right)
}=\sum_{l=0}^{k}A^{\left( k-l\right) }z^{\left( l\right) }=A^{\left(
0\right) }z^{\left( k\right) }+\sum_{l=0}^{k-1}A^{\left( k-l\right)
}z^{\left( l\right) }.
\end{equation*}%
The Adomian polynomials $z^{\left( k\right) }$ of $\ z=A^{-1}v$ are then
given by 
\begin{equation*}
A^{\left( 0\right) }z^{\left( k\right) }=v^{\left( k\right)
}-\sum_{l=0}^{k-1}A^{\left( k-l\right) }z^{\left( l\right) }.
\end{equation*}

\bigskip

Now to solve DAE system (\ref{M1}), let $Lp\left( t\right) =dp/dt$ and $%
L^{-2}p\left( t\right) =L^{-1}(L^{-1}p\left( t\right) ).$ Then $%
L^{-1}p\left( t\right) =\int_{0}^{t}p\left( t\right) dt$ and $L^{-2}p\left(
t\right) =\int_{0}^{t}\int_{0}^{t}p\left( t\right) dtdt.$

We solve the second equation of (\ref{M1}) for $dv/dt$ as 
\begin{equation}
\dfrac{dv}{dt}=M^{-1}\left( p\right) \Big(f\left( p,v\right) -G^{^{\text{%
\textsf{T}}}}\left( p\right) \lambda \Big)\text{.}  \label{M3-new}
\end{equation}%
Applying the operator $L^{-1}$ to both sides of the first equation of (\ref%
{M1}) and (\ref{M3-new}) then using initial conditions (\ref{M2}), we get 
\begin{equation}
\left\{ 
\begin{array}{l}
p=p_{0}+L^{-1}v\text{,} \\ 
v=v_{0}+L^{-1}\Big(M^{-1}\left( p\right) \Big(f\left( p,v\right) -G^{^{\text{%
\textsf{T}}}}\left( p\right) \lambda \Big)\Big).%
\end{array}%
\right.  \label{M4}
\end{equation}%
Then, we expand the solution components $p,v$ and $\lambda $ as 
\begin{equation}
p=\sum_{n=0}^{\infty }p^{\left( n\right) }\text{, \ }v=\sum_{n=0}^{\infty
}v^{\left( n\right) }\text{, \ }\lambda =\sum_{n=0}^{\infty }\lambda
^{\left( n\right) },  \label{M6}
\end{equation}%
where the unknowns $p^{\left( n\right) },v^{\left( n\right) }$ and $\lambda
^{\left( n\right) },$ $n=0,1,2,\ldots $ will be determined later by our
method. The nonlinear terms $f\left( p,v\right) $ and\ $g(p)$ are also
expanded in infinite series as 
\begin{equation}
f\left( p,v\right) =\sum_{n=0}^{\infty }f^{\left( n\right) }\text{, \ }%
g\left( p\right) =\sum_{n=0}^{\infty }g^{\left( n\right) }\text{ ,}
\label{M8}
\end{equation}%
where $f^{\left( n\right) }:=f^{\left( n\right) }\left( p^{\left( 0\right)
},v^{\left( 0\right) },\ldots ,p^{\left( n\right) },v^{\left( n\right)
}\right) $ and $g^{\left( n\right) }:=g^{\left( n\right) }\left( p^{\left(
0\right) },\ldots ,p^{\left( n\right) }\right) $ denote the Adomian
polynomials. Using (\ref{ADM-12}), the Adomian polynomials $f^{\left(
n\right) }$ and $g^{\left( n\right) },$ $n=0,1,2,\ldots $ can be written as%
\begin{equation}
f^{\left( n\right) }=\left\{ 
\begin{array}{l}
f\left( p^{\left( 0\right) },v^{\left( 0\right) }\right) ,\text{ \ }n=0, \\ 
\dfrac{1}{n}\sum_{i=1}^{n-1}i\left( \dfrac{\partial f^{\left( n-i\right) }}{%
\partial p^{\left( 0\right) }}\right) p^{\left( i\right) }+i\left( \dfrac{%
\partial f^{\left( n-i\right) }}{\partial v^{\left( 0\right) }}\right)
v^{\left( i\right) }+\left( \dfrac{\partial f^{\left( 0\right) }}{\partial
p^{\left( 0\right) }}\right) p^{\left( n\right) }+\left( \dfrac{\partial
f^{\left( 0\right) }}{\partial v^{\left( 0\right) }}\right) v^{\left(
n\right) },\text{ \ }n\geq 1,%
\end{array}%
\right.  \label{M10}
\end{equation}%
and%
\begin{equation}
g^{\left( n\right) }=\left\{ 
\begin{array}{l}
g\left( p^{\left( 0\right) }\right) ,\text{ \ }n=0, \\ 
\dfrac{1}{n}\sum_{i=1}^{n-1}i\left( \dfrac{\partial g^{\left( n-i\right) }}{%
\partial p^{\left( 0\right) }}\right) p^{\left( i\right) }+\left( \dfrac{%
\partial g^{\left( 0\right) }}{\partial p^{\left( 0\right) }}\right)
p^{\left( n\right) },\text{ \ }n\geq 1.%
\end{array}%
\right.  \label{M12}
\end{equation}%
Substituting expansions (\ref{M6}) into equations (\ref{M4}) and the third
equation of (\ref{M1}), we get%
\begin{equation}
\left\{ 
\begin{array}{l}
\sum_{n=0}^{\infty }p^{\left( n\right) }=p_{0}+\sum_{n=0}^{\infty
}L^{-1}v^{\left( n\right) }, \\ 
\sum_{n=0}^{\infty }v^{\left( n\right) }=v_{0}+\sum_{n=0}^{\infty }L^{-1}%
\Big(M^{-1}\left( p\right) \Big(f\left( u,v\right) -G^{^{\text{\textsf{T}}%
}}\left( p\right) \lambda \Big)\Big)^{\left( n\right) }, \\ 
0=\sum_{n=0}^{\infty }g^{\left( n\right) }.%
\end{array}%
\right.  \label{M14}
\end{equation}%
Choosing the initial terms $p^{\left( 0\right) }$ and $v^{\left( 0\right) }$
as 
\begin{equation}
p^{\left( 0\right) }=p_{0}\text{, \ }v^{\left( 0\right) }=v_{0},  \label{M18}
\end{equation}%
then comparing the general terms on the left and right hand sides of (\ref%
{M14}), we obtain the following recursion system 
\begin{equation}
\left\{ 
\begin{array}{l}
p^{\left( n\right) }=L^{-1}v^{\left( n-1\right) }, \\ 
v^{\left( n\right) }=L^{-1}\Big(M^{-1}\left( p\right) \Big(f\left(
p,v\right) -G^{^{\text{\textsf{T}}}}\left( p\right) \lambda \Big)\Big)%
^{\left( n-1\right) }, \\ 
0=g^{\left( n\right) },\text{ \ }n\geq 1.%
\end{array}%
\right.  \label{M20}
\end{equation}%
This system leads to the following recursion system 
\begin{equation}
\left\{ 
\begin{array}{l}
p^{\left( n\right) }=L^{-2}\Big(M^{-1}\left( p\right) \Big(f\left(
p,v\right) -G^{T}\left( p\right) \lambda \Big)\Big)^{\left( n-2\right) }, \\ 
0=g^{\left( n\right) },\text{ \ }n\geq 2,%
\end{array}%
\right.  \label{M22}
\end{equation}%
where%
\begin{equation}
v^{\left( n-1\right) }=Lp^{\left( n\right) }.  \label{M24}
\end{equation}%
Using expansions (\ref{M8}) and the previous proposition, we calculate the
right hand side of the first equation of system (\ref{M22}), and obtain the
following linear algebraic system for the unknowns $p^{\left( n\right) }$
and $L^{-2}\lambda ^{\left( n-2\right) }$ 
\begin{eqnarray}
M\left( p^{\left( 0\right) }\right) p^{\left( n\right) }+G^{^{\text{\textsf{T%
}}}}(p^{\left( 0\right) })L^{-2}\lambda ^{\left( n-2\right) } &=&R\text{,}
\label{M25} \\
G(p^{\left( 0\right) })p^{\left( n\right) } &=&S,\text{ \ }n\geq 2,  \notag
\end{eqnarray}%
where 
\begin{equation*}
R=L^{-2}\Bigg(f^{\left( n-2\right) }-\sum_{k=0}^{n-3}\Big(\left( G^{^{\text{%
\textsf{T}}}}(p)\right) ^{\left( n-2-k\right) }\lambda ^{\left( k\right)
}+\left( M\left( p\right) \right) ^{\left( n-2-k\right) }z^{\left( k\right) }%
\Big)\Bigg),
\end{equation*}%
and 
\begin{equation*}
S=-\dfrac{1}{n}\sum_{i=1}^{n-1}i\left( \dfrac{\partial g^{\left( n-i\right) }%
}{\partial p^{\left( 0\right) }}\right) p^{(i)}.
\end{equation*}%
The iterates $z^{\left( k\right) },$ $k=2,\ldots ,n-3$ are computed from%
\begin{equation}
M\left( p^{\left( 0\right) }\right) z^{\left( k\right) }=f^{\left( k\right)
}-G^{^{\text{\textsf{T}}}}(p^{\left( 0\right) })\lambda ^{\left( k\right)
}-\sum_{l=0}^{k-1}\Big(\left( G^{^{\text{\textsf{T}}}}(p)\right) ^{\left(
k-l\right) }\lambda ^{\left( l\right) }+\left( M\left( p\right) \right)
^{\left( k-l\right) }z^{\left( l\right) }\Big).  \label{M26}
\end{equation}%
In system (\ref{M25}),\ the right hand side depends only on previous
iterations $p^{\left( n-1\right) },\ldots ,p^{\left( 0\right) },\lambda
^{\left( n-3\right) },\ldots ,\lambda ^{\left( 0\right) }$ and $v^{\left(
n-1\right) },\ldots ,v^{\left( 0\right) }$. Since the Jacobian $G$ has full
row rank and the matrix$\ M$ is positive definite, then system (\ref{M25})
determines $p^{\left( n\right) }$ and $L^{-2}\lambda ^{\left( n-2\right) }$
uniquely for $n$ $\geq 2.$

One way to solve system (\ref{M25}) is to multiply the first equation of
this system from left by the matrix $G(p^{\left( 0\right) })M^{-1}(p^{\left(
0\right) }).$ Then substitute $G(p^{\left( 0\right) })p^{\left( n\right) }$
by its expression from the second equation of (\ref{M25}), to obtain the
following nonsingular algebraic system for the unknown $L^{-2}\lambda
^{\left( n-2\right) }$%
\begin{equation}
G(p^{\left( 0\right) })M^{-1}\left( p^{\left( 0\right) }\right) G^{^{\text{%
\textsf{T}}}}(p^{\left( 0\right) })L^{-2}\lambda ^{\left( n-2\right)
}=G(p^{\left( 0\right) })M^{-1}\left( p^{\left( 0\right) }\right) R-S.
\label{M27}
\end{equation}%
Since rank condition (\ref{Full-rank}) holds and $M$ is positive definite,
equation (\ref{M27}) can be solved uniquely for $L^{-2}\lambda ^{\left(
n-2\right) }$ to get%
\begin{equation}
L^{-2}\lambda ^{\left( n-2\right) }=\left( G(p^{\left( 0\right)
})M^{-1}\left( p^{\left( 0\right) }\right) G^{^{\text{\textsf{T}}%
}}(p^{\left( 0\right) })\right) ^{-1}\left( G(p^{\left( 0\right)
})M^{-1}\left( p^{\left( 0\right) }\right) R-S\right) .  \label{M29}
\end{equation}%
Now applying the operator $L^{2}$ to both sides of equation (\ref{M29}), we
can determine the unknown $\lambda ^{\left( n-2\right) }$ 
\begin{equation}
\lambda ^{\left( n-2\right) }=\left( G(p^{\left( 0\right) })M^{-1}\left(
p^{\left( 0\right) }\right) G^{^{\text{\textsf{T}}}}(p^{\left( 0\right)
})\right) ^{-1}\left( G(p^{\left( 0\right) })M^{-1}\left( p^{\left( 0\right)
}\right) L^{2}R-L^{2}S\right) .  \label{M30}
\end{equation}%
Then, substituting the expression of $L^{-2}\lambda ^{\left( n-2\right) }$
into the first equation of (\ref{M25}), we determine the unknown $p^{\left(
n\right) }$ 
\begin{equation}
p^{\left( n\right) }=M^{-1}\left( p^{\left( 0\right) }\right) R-M^{-1}\left(
p^{\left( 0\right) }\right) G^{^{\text{\textsf{T}}}}(p^{\left( 0\right)
})\left( G(p^{\left( 0\right) })M^{-1}\left( p^{\left( 0\right) }\right)
G^{^{\text{\textsf{T}}}}(p^{\left( 0\right) })\right) ^{-1}\left(
G(p^{\left( 0\right) })M^{-1}\left( p^{\left( 0\right) }\right) R-S\right) .
\label{M31}
\end{equation}%
Now, using equation (\ref{M24}) we can calculate $v^{\left( n-1\right) }.$
Finally, we obtain an approximate solution for DAE initial-value problem (%
\ref{M1})-(\ref{M2}) as 
\begin{equation}
p(t)=\displaystyle{{\sum_{n=0}^{n_{0}-1}}}p^{\left( n\right) }\text{, \ }%
v(t)=\displaystyle{{\sum_{n=0}^{n_{0}-2}}}v^{\left( n\right) }\text{, \ }%
\lambda (t)=\displaystyle{{\sum_{n=0}^{n_{0}-3}}}\lambda ^{\left( n\right) },
\label{MDAE-19}
\end{equation}%
where $n_{0}$ is the order of approximation of $p(t).$

\section{Application}

\label{Test problems}

In this section, we illustrate and demonstrate the effectiveness of our
technique to solve Euler-Lagrange equations (\ref{M1})-(\ref{M2}) which
describe the motion of constrained mechanical multibody systems. These
equations are known to be difficult to solve numerically because they are
index-three differential-algebraic equations (DAEs). Following the procedure
developed in the previous section, we first apply the Adomian decomposition
method (ADM) directly to these equations without using complex
transformations like index-reductions. Then, we expand the nonlinear terms
using the Adomian polynomials. Taking account of the index-three condition,
we derive a nonsingular linear algebraic recursion system for the expansion
components of the solution. Finally, by solving this algebraic system, we
obtain the solution of the Euler-Lagrange equations. As a test problem, we
consider the following example of constrained multibody system made up from
example 6.4 in \cite{Robot-example} which describes a two-link planar
robotic system, where the mass matrix is%
\begin{equation*}
M(\theta _{1},\theta _{2})=\left( 
\begin{array}{cc}
m_{1}l_{1}^{2}/3+m_{2}\left( l_{1}^{2}+l_{2}^{2}/3+l_{1}l_{2}\cos \theta
_{2}\right) & m_{2}\left( l_{2}^{2}/3+(1/2)l_{1}l_{2}\cos \theta _{2}\right)
\\ 
m_{2}\left( l_{2}^{2}/3+(1/2)l_{1}l_{2}\cos \theta _{2}\right) & 
m_{2}l_{2}^{2}/3%
\end{array}%
\right) .
\end{equation*}%
The force term is

\begin{equation*}
f\left( \theta _{1},\theta _{2},d\theta _{1}/dt,d\theta _{2}/dt\right)
=\left( 
\begin{array}{c}
\Big(l_{1}\cos \theta _{1}+l_{2}\cos \left( \theta _{1}+\theta _{2}\right) %
\Big)\left( d\theta _{1}/dt\right) -3\theta _{1} \\ 
\Big(l_{2}\cos \left( \theta _{1}+\theta _{2}\right) \Big)\left( d\theta
_{1}/dt\right) +\left( 1-\left( 3/2\right) \cos \theta _{2}\right) \theta
_{1}%
\end{array}%
\right) \text{, \ }
\end{equation*}%
and the constraint function is given by%
\begin{equation*}
g\left( \theta _{1},\theta _{2}\right) =l_{1}\sin \theta _{1}+l_{2}\sin
\left( \theta _{1}+\theta _{2}\right) .
\end{equation*}%
Taking $l_{1}=l_{2}=1,m_{1}=m_{2}=3$ and using the notation of (\ref{M1}),
with $n_{p}=2,$ $n_{\lambda }=1,$ $p=\left( p_{1},p_{2}\right) ^{^{\text{%
\textsf{T}}}}=\left( \theta _{1},\theta _{2}\right) ^{^{\text{\textsf{T}}}},$
$v=dp/dt$, we have 
\begin{equation*}
M\left( p\right) =\left( 
\begin{array}{cc}
5+3\cos p_{2} & \left( 3/2\right) \cos p_{2} \\ 
1+\left( 3/2\right) \cos p_{2} & 1%
\end{array}%
\right) ,
\end{equation*}%
and 
\begin{equation*}
f\left( p,v\right) =\left( 
\begin{array}{c}
\Big(\cos p_{1}+\cos \left( p_{1}+p_{2}\right) \Big)v_{1}-3p_{1} \\ 
\Big(\cos \left( p_{1}+p_{2}\right) \Big)v_{1}+\left( 1-\left( 3/2\right)
\cos p_{2}\right) p_{1}%
\end{array}%
\right) \text{. \ }
\end{equation*}%
The constraint function becomes 
\begin{equation*}
g\left( p\right) =\sin p_{1}+\sin \left( p_{1}+p_{2}\right) ,
\end{equation*}%
and its Jacobian $G\left( p\right) =$ $\left( \cos p_{1}+\cos \left(
p_{1}+p_{2}\right) ,\cos \left( p_{1}+p_{2}\right) \right) $ is full row
rank $n_{\lambda }=1.$

For the consistent initial conditions 
\begin{equation}
p\left( 0\right) =\left( 0,0\right) ^{^{\text{\textsf{T}}}}\text{, \ }%
v\left( 0\right) =\left( 1,-2\right) ^{^{\text{\textsf{T}}}},\text{\ }
\label{Ex-1-c}
\end{equation}%
the exact solution for this example is $p\left( t\right) =\left( \sin
t,-2\sin t\right) ^{^{\text{\textsf{T}}}}$, \ $v\left( t\right) =\left( \cos
t,-2\cos t\right) ^{^{\text{\textsf{T}}}}$ \ and $\lambda \left( t\right)
=\cos t.$

The Euler-Lagrange equations corresponding to this example form an
index-three DAE and therefore difficult to solve numerically. Using the
procedure developed in the previous section, system (\ref{M25}) can be
solved for $n=2,3,\ldots $ to reveal the dynamics of the mechanical system.

For $n=2,$ we get%
\begin{eqnarray}
p^{\left( 2\right) }+M^{-1}\left( p^{\left( 0\right) }\right) G^{^{\text{%
\textsf{T}}}}(p^{\left( 0\right) })L^{-2}\lambda ^{\left( 0\right) }
&=&L^{-2}M^{-1}\left( p^{\left( 0\right) }\right) f^{\left( 0\right) },
\label{Ex-6-a} \\
G(p^{\left( 0\right) })p^{\left( 2\right) } &=&0.  \label{Ex-6-b}
\end{eqnarray}%
Now, since 
\begin{equation*}
M^{-1}\left( p^{\left( 0\right) }\right) =\left( 2/7\right) \left( 
\begin{array}{cc}
2 & -5 \\ 
-5 & 16%
\end{array}%
\right) ,\ G(p^{\left( 0\right) })=\left( 2,1\right) ,
\end{equation*}%
\begin{equation*}
p^{\left( 0\right) }=\left( 
\begin{array}{c}
p_{1}^{\left( 0\right) } \\ 
p_{2}^{\left( 0\right) }%
\end{array}%
\right) =\left( 
\begin{array}{c}
0 \\ 
0%
\end{array}%
\right) ,\ v^{\left( 0\right) }=\left( 
\begin{array}{c}
v_{1}^{\left( 0\right) } \\ 
v_{2}^{\left( 0\right) }%
\end{array}%
\right) =\left( 
\begin{array}{c}
1 \\ 
-2%
\end{array}%
\right) ,
\end{equation*}%
and%
\begin{equation*}
f^{\left( 0\right) }=\left( 
\begin{array}{c}
\Big(\cos p_{1}^{\left( 0\right) }+\cos \left( p_{1}^{\left( 0\right)
}+p_{2}^{\left( 0\right) }\right) \Big)v_{1}^{\left( 0\right)
}-3p_{1}^{\left( 0\right) } \\ 
\Big(\cos \left( p_{1}^{\left( 0\right) }+p_{2}^{\left( 0\right) }\right) %
\Big)v_{1}^{\left( 0\right) }+\left( 1-\left( 3/2\right) \cos p_{2}^{\left(
0\right) }\right) p_{1}^{\left( 0\right) }%
\end{array}%
\right) =\left( 
\begin{array}{c}
2 \\ 
1%
\end{array}%
\right) ,
\end{equation*}%
equations (\ref{Ex-6-a})-(\ref{Ex-6-b}) reduce to%
\begin{eqnarray}
p^{\left( 2\right) }+\left( 2/7\right) \left( 
\begin{array}{cc}
2 & -5 \\ 
-5 & 16%
\end{array}%
\right) \left( 
\begin{array}{c}
2 \\ 
1%
\end{array}%
\right) L^{-2}\lambda ^{\left( 0\right) } &=&\left( 2/7\right) \left( 
\begin{array}{cc}
2 & -5 \\ 
-5 & 16%
\end{array}%
\right) \left( 
\begin{array}{c}
2 \\ 
1%
\end{array}%
\right) L^{-2}\left( 1\right) ,  \label{Ex-7-a} \\
\left( 2,1\right) p^{\left( 2\right) } &=&0.  \label{Ex-7-b}
\end{eqnarray}%
Multiplying system (\ref{Ex-7-a}) from left by the matrix $\left( 2,1\right) 
$ then using (\ref{Ex-7-b}), we get%
\begin{equation}
\left( 2/7\right) \left( 2,1\right) \left( 
\begin{array}{cc}
2 & -5 \\ 
-5 & 16%
\end{array}%
\right) \left( 
\begin{array}{c}
2 \\ 
1%
\end{array}%
\right) L^{-2}\lambda ^{\left( 0\right) }=\left( 2/7\right) \left(
2,1\right) \left( 
\begin{array}{cc}
2 & -5 \\ 
-5 & 16%
\end{array}%
\right) \left( 
\begin{array}{c}
2 \\ 
1%
\end{array}%
\right) L^{-2}\left( 1\right) ,  \label{Ex-8}
\end{equation}%
which leads

\begin{equation}
\lambda ^{\left( 0\right) }=1.  \label{Ex-9}
\end{equation}%
Now substituting the value of $\lambda ^{\left( 0\right) }$ from (\ref{Ex-9}%
) into (\ref{Ex-7-a}), we obtain%
\begin{equation}
p^{\left( 2\right) }+\left( 2/7\right) \left( 
\begin{array}{cc}
2 & -5 \\ 
-5 & 16%
\end{array}%
\right) \left( 
\begin{array}{c}
2 \\ 
1%
\end{array}%
\right) L^{-2}\left( 1\right) =\left( 2/7\right) \left( 
\begin{array}{cc}
2 & -5 \\ 
-5 & 16%
\end{array}%
\right) \left( 
\begin{array}{c}
2 \\ 
1%
\end{array}%
\right) L^{-2}\left( 1\right) ,  \label{Ex-10}
\end{equation}%
which gives 
\begin{equation}
p^{\left( 2\right) }=\left( 
\begin{array}{c}
0 \\ 
0%
\end{array}%
\right) ,  \label{Ex-11}
\end{equation}%
and using (\ref{M24}), we have 
\begin{equation}
v^{\left( 1\right) }=Lp^{\left( 2\right) }=\left( 
\begin{array}{c}
0 \\ 
0%
\end{array}%
\right) .  \label{Ex-12}
\end{equation}%
For $n=3$, we have%
\begin{eqnarray}
p^{\left( 3\right) }+M^{-1}\left( p^{\left( 0\right) }\right) G^{^{\text{%
\textsf{T}}}}(p^{\left( 0\right) })L^{-2}\lambda ^{\left( 1\right) }
&=&L^{-2}M^{-1}\left( p^{\left( 0\right) }\right) \Bigg(f^{\left( 1\right) }
\notag \\
&&-\left( G^{^{\text{\textsf{T}}}}(p)\right) ^{\left( 1\right) }\lambda
^{\left( 0\right) }-\left( M\left( p\right) \right) ^{\left( 1\right)
}z^{\left( 0\right) }\Bigg)\text{,}  \label{Ex-13} \\
G(p^{\left( 0\right) })p^{\left( 3\right) } &=&0.  \notag
\end{eqnarray}%
Now, since 
\begin{equation*}
\left( M\left( p\right) \right) ^{\left( 1\right) }=\left( 
\begin{array}{cc}
-3p_{2}^{\left( 1\right) }\sin p_{2}^{\left( 0\right) } & -\left( 3/2\right)
p_{2}^{\left( 1\right) }\sin p_{2}^{\left( 0\right) } \\ 
-\left( 3/2\right) p_{2}^{\left( 1\right) }\sin p_{2}^{\left( 0\right) } & 0%
\end{array}%
\right) =0,\ \left( G^{^{\text{\textsf{T}}}}(p)\right) ^{\left( 1\right)
}=0,\ 
\end{equation*}%
and%
\begin{equation*}
f^{\left( 1\right) }=\left( 
\begin{array}{c}
-3 \\ 
-1/2%
\end{array}%
\right) t,
\end{equation*}%
system (\ref{Ex-13}) reduces to%
\begin{eqnarray}
p^{\left( 3\right) }+\left( 2/7\right) \left( 
\begin{array}{cc}
2 & -5 \\ 
-5 & 16%
\end{array}%
\right) \left( 
\begin{array}{c}
2 \\ 
1%
\end{array}%
\right) L^{-2}\lambda ^{\left( 1\right) } &=&\left( 2/7\right) \left( 
\begin{array}{cc}
2 & -5 \\ 
-5 & 16%
\end{array}%
\right) \left( 
\begin{array}{c}
-3 \\ 
-1/2%
\end{array}%
\right) L^{-2}\left( t\right) ,  \label{Ex-14-a} \\
&&\left( 2,1\right) p^{\left( 3\right) }=0.  \label{Ex-14-b}
\end{eqnarray}%
Multiplying system (\ref{Ex-14-a}) from left by the matrix $(2,1)$ then
using (\ref{Ex-14-b}), we get%
\begin{eqnarray}
\left( 2/7\right) \left( 2,1\right) \left( 
\begin{array}{cc}
2 & -5 \\ 
-5 & 16%
\end{array}%
\right) \left( 
\begin{array}{c}
2 \\ 
1%
\end{array}%
\right) L^{-2}\lambda ^{\left( 1\right) } &=&\left( 2/7\right) \left(
2,1\right) \left( 
\begin{array}{cc}
2 & -5 \\ 
-5 & 16%
\end{array}%
\right) \left( 
\begin{array}{c}
-3 \\ 
-1/2%
\end{array}%
\right) L^{-2}\left( t\right) ,  \notag \\
&=&0,  \label{Ex-15}
\end{eqnarray}%
which gives

\begin{equation}
\lambda ^{\left( 1\right) }=0.  \label{Ex-17}
\end{equation}%
Now substituting the value of $\lambda ^{\left( 1\right) }$ from (\ref{Ex-17}%
) into (\ref{Ex-14-a}), we obtain%
\begin{equation}
p^{\left( 3\right) }=\left( 2/7\right) \left( 
\begin{array}{cc}
2 & -5 \\ 
-5 & 16%
\end{array}%
\right) \left( 
\begin{array}{c}
-3 \\ 
-1/2%
\end{array}%
\right) L^{-2}\left( t\right) ,  \label{Ex-18}
\end{equation}%
which gives 
\begin{equation}
p^{\left( 3\right) }=-\dfrac{t^{3}}{3!}\left( 
\begin{array}{c}
1 \\ 
-2%
\end{array}%
\right) ,  \label{Ex-19}
\end{equation}%
and 
\begin{equation}
v^{\left( 2\right) }=Lp^{\left( 3\right) }=-\dfrac{t^{2}}{2!}\left( 
\begin{array}{c}
1 \\ 
-2%
\end{array}%
\right) .  \label{Ex-20}
\end{equation}%
For $n=4$, we have%
\begin{eqnarray}
p^{\left( 4\right) }+M^{-1}\left( p^{\left( 0\right) }\right) G^{^{\text{%
\textsf{T}}}}(p^{\left( 0\right) })L^{-2}\lambda ^{\left( 2\right) }
&=&L^{-2}M^{-1}\left( p^{\left( 0\right) }\right) \Bigg(f^{\left( 2\right) }
\notag \\
&&-\left( G^{^{\text{\textsf{T}}}}(p)\right) ^{\left( 2\right) }\lambda
^{\left( 0\right) }-\left( M\left( p\right) \right) ^{\left( 2\right)
}z^{\left( 0\right) }  \notag \\
&&-\left( G^{^{\text{\textsf{T}}}}(p)\right) ^{\left( 1\right) }\lambda
^{\left( 1\right) }-\left( M\left( p\right) \right) ^{\left( 1\right)
}z^{\left( 1\right) }\Bigg)\text{,}  \label{Ex-21} \\
G(p^{\left( 0\right) })p^{\left( 4\right) } &=&-\dfrac{1}{4}%
\sum_{i=1}^{3}i\left( \dfrac{\partial g^{\left( 4-i\right) }}{\partial
p^{\left( 0\right) }}\right) p^{(i)},  \notag
\end{eqnarray}%
and 
\begin{equation*}
f^{\left( 2\right) }=\left( 
\begin{array}{c}
-2 \\ 
-1%
\end{array}%
\right) t^{2},\ \left( G^{^{\text{\textsf{T}}}}(p)\right) ^{\left( 2\right)
}=\left( 
\begin{array}{c}
-1 \\ 
-1/2%
\end{array}%
\right) t^{2}.
\end{equation*}%
The iterates $z^{\left( 0\right) }$ is calculated from%
\begin{equation*}
M\left( p^{\left( 0\right) }\right) z^{\left( 0\right) }=f^{\left( 0\right)
}-G^{^{\text{\textsf{T}}}}(p^{\left( 0\right) })\lambda ^{\left( 0\right)
}=\left( 
\begin{array}{c}
2 \\ 
1%
\end{array}%
\right) -\left( 
\begin{array}{c}
2 \\ 
1%
\end{array}%
\right) \left( 1\right) =0,
\end{equation*}%
which gives $z^{\left( 0\right) }=\left( 
\begin{array}{c}
0 \\ 
0%
\end{array}%
\right) .$ Thus (\ref{Ex-21}) becomes 
\begin{eqnarray}
p^{\left( 4\right) }+\left( 2/7\right) \left( 
\begin{array}{cc}
2 & -5 \\ 
-5 & 16%
\end{array}%
\right) \left( 
\begin{array}{c}
2 \\ 
1%
\end{array}%
\right) L^{-2}\lambda ^{\left( 2\right) } &=&\left( 2/7\right) \left( 
\begin{array}{cc}
2 & -5 \\ 
-5 & 16%
\end{array}%
\right) \left( 
\begin{array}{c}
2 \\ 
1%
\end{array}%
\right) L^{-2}\left( -1/2t^{2}\right) ,  \label{Ex-22-a} \\
(2,1)p^{\left( 4\right) } &=&0.  \label{Ex-22-b}
\end{eqnarray}%
Multiplying system \ (\ref{Ex-22-a}) from left by the matrix $(2,1)$ then
using (\ref{Ex-22-b}), we get%
\begin{equation}
\left( 2/7\right) \left( 2,1\right) \left( 
\begin{array}{cc}
2 & -5 \\ 
-5 & 16%
\end{array}%
\right) \left( 
\begin{array}{c}
2 \\ 
1%
\end{array}%
\right) L^{-2}\lambda ^{\left( 2\right) }=\left( 2/7\right) \left(
2,1\right) \left( 
\begin{array}{cc}
2 & -5 \\ 
-5 & 16%
\end{array}%
\right) \left( 
\begin{array}{c}
2 \\ 
1%
\end{array}%
\right) L^{-2}\left( -1/2t^{2}\right) ,  \label{Ex-24}
\end{equation}%
which gives 
\begin{equation}
\lambda ^{\left( 2\right) }=-1/2t^{2}.  \label{Ex-25}
\end{equation}%
Now substituting the value of $\lambda ^{\left( 2\right) }$ from (\ref{Ex-25}%
) into (\ref{Ex-22-a}), we obtain%
\begin{equation*}
p^{\left( 4\right) }+\left( 2/7\right) \left( 
\begin{array}{cc}
2 & -5 \\ 
-5 & 16%
\end{array}%
\right) \left( 
\begin{array}{c}
2 \\ 
1%
\end{array}%
\right) L^{-2}\left( -1/2t^{2}\right) =\left( 2/7\right) \left( 
\begin{array}{cc}
2 & -5 \\ 
-5 & 16%
\end{array}%
\right) \left( 
\begin{array}{c}
2 \\ 
1%
\end{array}%
\right) L^{-2}\left( -1/2t^{2}\right) ,
\end{equation*}%
which gives%
\begin{equation}
p^{\left( 4\right) }=\left( 
\begin{array}{c}
0 \\ 
0%
\end{array}%
\right) ,  \label{Ex-26}
\end{equation}%
and using (\ref{M24}), we have 
\begin{equation}
v^{\left( 3\right) }=Lp^{\left( 4\right) }=\left( 
\begin{array}{c}
0 \\ 
0%
\end{array}%
\right) .  \label{Ex-27}
\end{equation}%
For $n=5$, we have%
\begin{eqnarray}
p^{\left( 5\right) }+M^{-1}\left( p^{\left( 0\right) }\right) G^{^{\text{%
\textsf{T}}}}(p^{\left( 0\right) })L^{-2}\lambda ^{\left( 3\right) }
&=&L^{-2}M^{-1}\left( p^{\left( 0\right) }\right) \Bigg(f^{\left( 3\right) }
\notag \\
&&-\sum_{k=0}^{2}\Big(\left( G^{^{\text{\textsf{T}}}}(p)\right) ^{\left(
3-k\right) }\lambda ^{\left( k\right) }+\left( M\left( p\right) \right)
^{\left( 3-k\right) }z^{\left( k\right) }\Big)\Bigg)\text{,}  \label{Ex-28}
\\
G(p^{\left( 0\right) })p^{\left( 5\right) } &=&-\dfrac{1}{5}%
\sum_{i=1}^{4}i\left( \dfrac{\partial g^{\left( 5-i\right) }}{\partial
p^{\left( 0\right) }}\right) p^{(i)},  \notag
\end{eqnarray}%
where 
\begin{equation*}
f^{\left( 3\right) }=\left( 
\begin{array}{c}
1/2 \\ 
37/12%
\end{array}%
\right) t^{3},\ \left( G^{^{\text{\textsf{T}}}}(p)\right) ^{\left( 3\right)
}=0.
\end{equation*}%
\begin{equation*}
\left( M(p)\right) ^{\left( 2\right) }=\left( 
\begin{array}{cc}
-6t^{2} & -3t^{2} \\ 
-3t^{2} & 0%
\end{array}%
\right) .
\end{equation*}%
The term $z^{\left( 1\right) }$ is calculated from

\begin{eqnarray*}
M\left( p^{\left( 0\right) }\right) z^{\left( 1\right) } &=&f^{\left(
1\right) }-\left( G^{^{\text{\textsf{T}}}}(p)\right) ^{\left( 1\right)
}\lambda ^{\left( 0\right) }-G^{^{\text{\textsf{T}}}}(p^{\left( 0\right)
})\lambda ^{\left( 1\right) }-\left( M\left( p\right) \right) ^{\left(
1\right) }z^{\left( 0\right) } \\
&=&f^{\left( 1\right) }=\left( 
\begin{array}{c}
-3 \\ 
-1/2%
\end{array}%
\right) t,
\end{eqnarray*}%
which gives $z^{\left( 1\right) }=\left( 
\begin{array}{c}
-1 \\ 
2%
\end{array}%
\right) t.$ Thus system (\ref{Ex-28}) becomes 
\begin{eqnarray}
p^{\left( 5\right) }+\left( 2/7\right) \left( 
\begin{array}{cc}
2 & -5 \\ 
-5 & 16%
\end{array}%
\right) \left( 
\begin{array}{c}
2 \\ 
1%
\end{array}%
\right) L^{-2}\lambda ^{\left( 3\right) } &=&\left( 2/7\right) \left( 
\begin{array}{cc}
2 & -5 \\ 
-5 & 16%
\end{array}%
\right) \left( 
\begin{array}{c}
1/2 \\ 
1/12%
\end{array}%
\right) L^{-2}\left( t^{3}\right) ,  \label{Ex-29-a} \\
(2,1)p^{\left( 5\right) } &=&0.  \label{Ex-29-b}
\end{eqnarray}%
Multiplying system (\ref{Ex-29-a}) from left by the matrix $(2,1)$ then
using (\ref{Ex-29-b}), we get%
\begin{eqnarray}
\left( 2/7\right) \left( 2,1\right) \left( 
\begin{array}{cc}
2 & -5 \\ 
-5 & 16%
\end{array}%
\right) \left( 
\begin{array}{c}
2 \\ 
1%
\end{array}%
\right) L^{-2}\lambda ^{\left( 3\right) } &=&\left( 2/7\right) \left(
2,1\right) \left( 
\begin{array}{cc}
2 & -5 \\ 
-5 & 16%
\end{array}%
\right) \left( 
\begin{array}{c}
1/2 \\ 
1/12%
\end{array}%
\right) L^{-2}\left( t^{3}\right) ,  \notag \\
&=&0,  \label{Ex-30}
\end{eqnarray}%
which gives 
\begin{equation}
\lambda ^{\left( 3\right) }=0.  \label{Ex-31}
\end{equation}%
Now substituting the value of $\lambda ^{\left( 3\right) }$ from (\ref{Ex-31}%
) into (\ref{Ex-29-a}), we obtain%
\begin{equation}
p^{\left( 5\right) }=\left( 2/7\right) \left( 
\begin{array}{cc}
2 & -5 \\ 
-5 & 16%
\end{array}%
\right) \left( 
\begin{array}{c}
1/2 \\ 
1/12%
\end{array}%
\right) L^{-2}\left( t^{3}\right) ,  \label{Ex-32}
\end{equation}%
which gives%
\begin{equation}
p^{\left( 5\right) }=\dfrac{t^{5}}{5!}\left( 
\begin{array}{c}
1 \\ 
-2%
\end{array}%
\right) ,  \label{Ex-33}
\end{equation}%
and using (\ref{M24}), we have 
\begin{equation}
v^{\left( 4\right) }=Lp^{\left( 5\right) }=\dfrac{t^{4}}{4!}\left( 
\begin{array}{c}
1 \\ 
-2%
\end{array}%
\right) .  \label{Ex-34}
\end{equation}%
Continuing this process until $n=8$, we obtain the following ADM solution 
\begin{equation}
p\left( t\right) =t\left( 
\begin{array}{c}
1 \\ 
-2%
\end{array}%
\right) -\dfrac{t^{3}}{3!}\left( 
\begin{array}{c}
1 \\ 
-2%
\end{array}%
\right) +\dfrac{t^{5}}{5!}\left( 
\begin{array}{c}
1 \\ 
-2%
\end{array}%
\right) -\dfrac{t^{7}}{7!}\left( 
\begin{array}{c}
1 \\ 
-2%
\end{array}%
\right) ,  \label{Ex-35}
\end{equation}%
\begin{equation}
v\left( t\right) =\left( 
\begin{array}{c}
1 \\ 
-2%
\end{array}%
\right) -\dfrac{t^{2}}{2!}\left( 
\begin{array}{c}
1 \\ 
-2%
\end{array}%
\right) +\dfrac{t^{4}}{4!}\left( 
\begin{array}{c}
1 \\ 
-2%
\end{array}%
\right) -\dfrac{t^{6}}{6!}\left( 
\begin{array}{c}
1 \\ 
-2%
\end{array}%
\right) ,  \label{EX-35a}
\end{equation}%
and 
\begin{equation}
\lambda \left( t\right) =1-\dfrac{t^{2}}{2!}+\dfrac{t^{4}}{4!}-\dfrac{t^{6}}{%
6!}.  \label{Ex-36}
\end{equation}%
These are the first few terms of Taylor expansions, around $t=0,$ of 
\begin{equation}
p\left( t\right) =\left( \sin t,-2\sin t\right) ^{^{\text{\textsf{T}}}}\text{%
, \ \ }v\left( t\right) =\left( \cos t,-2\cos t\right) ^{^{\text{\textsf{T}}%
}}\text{, \ \ }\lambda \left( t\right) =\cos t,  \label{Ex-37}
\end{equation}%
which is the exact solution of the problem in this example.

\section{Discussion}

\label{Discussion}The Euler-Lagrange equations are known to be difficult to
solve numerically. The reason is that they form an index-three system of
differential-algebraic equations (DAEs). In this paper, we propose a novel
technique that applies the Adomian decomposition method (ADM) directly to
solve the Euler-Lagrange equations. This technique has successfully handled
these equations without the need for complex transformations like
index-reductions. This method transforms these equations into easily
solvable algebraic systems for the expansion components of the solution. To
illustrate the effectiveness of the proposed technique, an example of the
Euler-Lagrange equations describing a two-link planar robot system is
solved. This example shows that the ADM is a simple powerful tool to obtain
the exact or approximate solutions of the Euler-Lagrange equations.

\section{Conclusion}

\label{Conclusion}

This work presents the analytical solution of the Euler-Lagrange equations
using the ADM. A \ procedure for solving these is presented. The technique
was tested on an example of the Euler-Lagrange equations that describes a
two-link robot system. The results obtained show that the proposed method
can be applied to solve the Euler-Lagrange equations efficiently to obtain
the exact or an approximate solution. On the one hand, it is important to
note that these types of equations are difficult to solve and on the other,
the direct application of the ADM was able to solve the Euler-Lagrange
equations. Also, it is important to note that, our technique does not make
transformations to the equations before applying the ADM to them. The
technique is based on a straightforward procedure that can be programmed in
Maple or Mathematica to simulate real application problems. Finally, further
work is needed to apply a multistage ADM form to solve the Euler-Lagrange
equations and other semi-explicit nonlinear higher-index DAEs. \ 

\textbf{Conflict of Interests}

The author declares that there is no conflict of interests regarding the
publication of this paper.

\end{document}